\documentclass[10pt]{article}
\usepackage[utf8]{inputenc}
\usepackage[T1]{fontenc}
\usepackage{amsmath}
\usepackage{amsfonts}
\usepackage{amssymb}
\usepackage[version=4]{mhchem}
\usepackage{stmaryrd}
\usepackage{graphicx}
\usepackage[export]{adjustbox}
\graphicspath{ {./images/} }
\usepackage{hyperref}
\hypersetup{colorlinks=true, linkcolor=blue, filecolor=magenta, urlcolor=cyan,}
\title{Accelerating the Hypergeometric Function with the Beta Integral to Derive New Infinite Series for $\pi$ and Values of the Gamma Function}

\author{Cetin Hakimoglu-Brown\\
mathemails@protonmail.com}
  
\newcommand\longdiv[1]{\overline{\smash{)}#1}}

\begin{document}
\maketitle


\begin{abstract}
The beta integral is applied to accelerate the hypergeometric function $2 F 1\left\{1, B; C ; w\right\}$ to derive new infinite series for  constants such as $\pi$ and values of the gamma function. A compendium of new infinite series is given. Ramanujan-like formulas for pi are also derived based on elementary inverse trigonometric functions, including a formula with rational values that  adds 2.5 digits per terms, which makes the series much more compact than similar formulas in the existing literature.
\end{abstract}

\section{Introduction}
\numberwithin{equation}{section}

At around the  '90s and early 2000s, David \& Peter Borwein et al. [2] Travis Sherman [3],  
K. A. Driver et al. [13], Batir [14]  Christian Krattenthaler et al. [4], and Boris Gour'evitch \& Jes'us Guillera Goyanes [1], Bellard [6] derived various binomial-like infinite series for constants such as $\pi$, by using the beta integral as applied to inverse trigonometric integrals. This paper continues on such earlier work, and extends it to the Euler integral of the hypergeometric function, in order to make the the algorithm more comprehensive and generalized.  

Examples of formulas that can be derived though the algorithm include  the so-called '324 formula' for the exceptionally fast rate of convergence:

\begin{equation}
\pi=\frac{\sqrt{3}}{60} \sum_{n=0}^{\infty} \frac{(2 n) !(130 n+109)}{\left(\frac{7}{6}\right)_{n}\left(\frac{11}{6}\right)_{n}(-1296)^{n}}
\end{equation}

And

\begin{equation}
\frac{7 \cdot 3^{1 / 2}(\Gamma(1 / 3))^{3}}{\pi \cdot 2^{1 / 3}}=\sum_{n=0}^{\infty} \frac{\left(\frac{2}{3}\right)_{n}\left(\frac{1}{2}\right)_{n}(102 n+59)}{\left(\frac{13}{12}\right)_{n}\left(\frac{19}{12}\right)_{n}(-288)^{n}}
\end{equation}

These formulas were derived with the following integrals which we shall prove in this paper:

\begin{equation}
\frac{z \cdot \Gamma(a+b+2)}{\Gamma(a+1) \Gamma(b+1)} \int_{0}^{1} \frac{x^{a}(1-x)^{b}}{z-x^{k}(1-x)^{s}} d x=\sum_{n=0}^{\infty} \frac{(a+1)_{k n}(b+1)_{s n}}{(a+b+2)_{(k+s) n}(z)^{n}}
\end{equation}

We'll also prove the more comprehensive integral-polynomial relation. This is the so called 'machinery' behind these infinite series:

\begin{equation}
\begin{aligned}
& \frac{\Gamma(a+1) \Gamma(b+1)}{Z \cdot \Gamma(a+b+2)} \sum_{n=0}^{\infty}\left(\frac{(a+1)_{k n}(b+1)_{s n}}{(a+b+2)_{(k+s) n}(z)^{n}}(w)\right)= \\
& \int_{0}^{1} \frac{Q(x) x^{a}(1-x)^{b}}{z-x^{k}(1-x)^{s}} d x=\int_{0}^{1} \frac{x^{a}(1-x)^{b}}{P(x)} d x=c \\
& w=\left\{\begin{array}{c}
a_{p} \prod_{g=1}^{p}\left(\frac{a+g+k n}{a+b+g+1+(k+s) n}\right)+a_{p-1} \prod_{g=1}^{p-1}\left(\frac{a+g+k n}{a+b+g+1+(k+s) n}\right)+ \\
\ldots+a_{1}\left(\frac{a+1+k n}{a+b+2+(k+s) n}\right)+a_{0}
\end{array}\right\}
\end{aligned}
\end{equation}

This paper was originally composed in 2009 and published in its current form to SSRN in 2021 [8]. The final section includes new results. Later similar results involving binomials since the publication of this paper include John M. Campbell et al., 2019 [9], Wenchang Chu, 2011 [10],  Kunle Adegok  et al., 2023 [11]. The '324 formula' and the exact implementation of the algorithm used here was hitherto unknown, and independently re-derived by Zhi-Wei Sun, 2023 [12]. 

\section{Derivation of Algorithm}
\numberwithin{equation}{section}
\addtocounter{equation}{-1}

We will begin by deriving formula (1.3) and apply it to derive (1.4).

Consider the beta function:

$$
\int_{0}^{1} x^{k n+a}(1-x)^{s n+b} d x=\frac{\Gamma(k n+a+1) \Gamma(s n+b+1)}{\Gamma((k+s) n+a+b+2)}
$$

Then divide both sides by $z^{n}$ and take the summation

$$
\sum_{n=0}^{\infty}\left(\int_{0}^{1} \frac{x^{k n+a}(1-x)^{s n+b}}{z^{n}} d x\right)=\sum_{n=0}^{\infty} \frac{\Gamma(k n+a+1) \Gamma(s n+b+1)}{\Gamma(a+b+2+n(s+k)) \cdot z^{n}}
$$

This is equal to:

\begin{equation}
\sum_{n=0}^{\infty}\left(\int_{0}^{1} x^{a}(1-x)^{b}\left(\frac{x^{k}(1-x)^{s}}{z}\right)^{n} d x\right)=\sum_{n=0}^{\infty} \frac{\Gamma(k n+a+1) \Gamma(s n+b+1)}{\Gamma(a+b+2+n(s+k)) \cdot z^{n}}
\end{equation}

Summing an infinite geometric series:

\begin{equation}
\sum_{n=0}^{\infty}\left(\int_{0}^{1} x^{a}(1-x)^{b}\left(\frac{x^{k}(1-x)^{s}}{z}\right)^{n} d x\right)=z \int_{0}^{1} \frac{x^{a}(1-x)^{b}}{z-x^{k}(1-x)^{s}} d x
\end{equation}

Using the identity: $(\Gamma(w))(w)_{n}=\Gamma(w+n)$ we derive:

\begin{equation}
\frac{\Gamma(k n+a+1) \Gamma(s n+b+1)}{\Gamma((k+s) n+a+b+2)}=\frac{(a+1)_{k n}(b+1)_{s n} \Gamma(a+1) \Gamma(b+1)}{(a+b+2)_{(k+s) n} \Gamma(a+b+2)}
\end{equation}

Plugging equations (2.1) and (2.2) into (2.0) and re-arranging terms we get equation (1.3)

To make the integral useful for computing constants begin with a so called 'seed integral'. Such an integral must be in the form (for linear $P(x)$ is the Euler integral for the hypergeoemtric function):

$$
\int_{0}^{1} \frac{x^{a}(1-x)^{b}}{P(x)} d x=c \text { where } P(x) \text { is a polynomial and } c \text { is a constant }
$$

Then choose integer values of $k$ and $s$ such that the polynomial expansion of $z-x^{k}(1-x)^{s}$ is divisible by $P(x)$. The value of $z$ is determined by computing the quotient polynomial $Q(x)$ via polynomial long division such that the following identity is obtained:

$$
\int_{0}^{1} \frac{x^{a}(1-x)^{b}}{P(x)} d x=\int_{0}^{1} \frac{Q(x) x^{a}(1-x)^{b}}{z-x^{k}(1-x)^{s}} d x
$$

$Q(x)$ can be generalized to: $Q(x)=a_{p} x^{p}+a_{p-1} x^{p-1} \ldots a_{1} x+a_{0}$

Multiplying $\int_{0}^{1} \frac{x^{a}(1-x)^{b}}{z-x^{k}(1-x)^{s}} d x$ by $Q(x)$ and splitting the $a_{p}$ coefficients of $Q(x)$ yields:

$$
\left\{\begin{array}{l}
a_{p} \int_{0}^{1} \frac{x^{a+p}(1-x)^{b}}{z-x^{k}(1-x)^{s}} d x+a_{p-1} \int_{0}^{1} \frac{x^{a+p-1}(1-x)^{b}}{z-x^{k}(1-x)^{s}} d x \\
\ldots a_{1} \int_{0}^{1} \frac{x^{a+1}(1-x)^{b}}{z-x^{k}(1-x)^{s}} d x+a_{0} \int_{0}^{1} \frac{x^{a}(1-x)^{b}}{z-x^{k}(1-x)^{s}} d x
\end{array}\right\}=\int_{0}^{1} \frac{x^{a}(1-x)^{b}}{P(x)} d x
$$

Each of these $a_{p} \ldots a_{0}$ split terms can be plugged into equation (1.3) :

\begin{equation}
\left\{\begin{array}{l}
a_{p} \int_{0}^{1} \frac{x^{a+p}(1-x)^{b}}{z-x^{k}(1-x)^{s}} d x=\frac{a_{p} \cdot \Gamma(a+p+1) \Gamma(b+1)}{z \cdot \Gamma(a+b+p+2)} \sum_{n=0}^{\infty} \frac{(a+p+1)_{k n}(b+1)_{s n}}{(a+b+p+2)_{(k+s) n}(z)^{n}} \\
a_{p-1} \int_{0}^{1} \frac{x^{a+p-1}(1-x)^{b}}{z-x^{k}(1-x)^{s}} d x=\frac{a_{p-1} \cdot \Gamma(a+p) \Gamma(b+1)}{z \cdot \Gamma(a+b+p+1)} \sum_{n=0}^{\infty} \frac{(a+p)_{k n}(b+1)_{s n}}{(a+b+p+1)_{(k+s) n}(z)^{n}} \\
\text { - } \\
a_{1} \int_{0}^{1} \frac{x^{a+1}(1-x)^{b}}{z-x^{k}(1-x)^{s}} d x=\frac{a_{1} \cdot \Gamma(a+2) \Gamma(b+1)}{z \cdot \Gamma(a+b+3)} \sum_{n=0}^{\infty} \frac{(a+2)_{k n}(b+1)_{s n}}{(a+b+3)_{(k+s) n}(z)^{n}} \\
a_{0} \int_{0}^{1} \frac{x^{a}(1-x)^{b}}{z-x^{k}(1-x)^{s}} d x=\frac{a_{0} \cdot \Gamma(a+1) \Gamma(b+1)}{z \cdot \Gamma(a+b+2)} \sum_{n=0}^{\infty} \frac{(a+1)_{k n}(b+1)_{s n}}{(a+b+2)_{(k+s) n}(z)^{n}}
\end{array}\right\}
\end{equation}

Also note the following pochhammer and gamma function identities:

\begin{equation}
\begin{aligned}
& (a+b+p+2)_{(k+s) n}=\left(\prod_{g=1}^{p}\left(\frac{a+b+g+1+(k+s) n}{a+b+g+1}\right)\right)(a+b+2)_{(k+s) n} \\
& (a+p+1)_{k n}=\left(\prod_{g=1}^{p}\left(\frac{a+g+k n}{a+g}\right)\right)(a+1)_{(k n)} \\
& \Gamma(a+b+p+2)=\left(\prod_{g=1}^{p}(a+b+g+1)\right) \Gamma(a+b+2) \\
& \Gamma(a+p+1)=\left(\prod_{g=1}^{p}(a+g)\right) \Gamma(a+1)
\end{aligned}
\end{equation}

Putting it all together using (2.3) and (2.4) gives:

$$
\begin{aligned}
& \frac{a_{p} \cdot \Gamma(a+p+1) \Gamma(b+1)}{z \cdot \Gamma(a+b+p+2)} \sum_{n=0}^{\infty} \frac{(a+p+1)_{k n}(b+1)_{s n}}{(a+b+p+2)_{(k+s) n}(z)^{n}}= \\
& \frac{a_{p}\left(\prod_{g=1}^{p}(a+g)\right) \Gamma(a+1) \Gamma(b+1)}{z \cdot\left(\prod_{g=1}^{p}(a+b+g+1)\right) \Gamma(a+b+2)} \sum_{n=0}^{\infty}\left(\frac{\left(\prod_{g=1}^{p}\left(\frac{a+g+k n}{a+g}\right)\right)(a+1)_{k n}(b+1)_{s n}}{\left.\left(\prod_{g=1}^{p}\left(\frac{a+b+g+1+(k+s) n}{a+b+g+1}\right)\right)(a+b+2)_{(k+s) n}(z)^{n}\right)}\right)\end{aligned}
$$

Also note how the $\prod_{g=1}^{p}(a+b+g+1)$ and $\prod_{g=1}^{p}(a+g)$ are eliminated such that we simplify and obtain:

\begin{equation}
\begin{aligned}
& \frac{a_{p} \Gamma(a+1) \Gamma(b+1)}{z \cdot \Gamma(a+b+2)} \sum_{n=0}^{\infty}\left(\frac{\left(\prod_{g=1}^{p}(a+g+k n)\right)(a+1)_{k n}(b+1)_{s n}}{\left(\prod_{g=1}^{p}(a+b+g+1+(k+s) n)\right)(a+b+2)_{(k+s) n}(z)^{n}}\right)= \\
& a_{p} \int_{0}^{1} \frac{x^{a+p}(1-x)^{b}}{z-x^{k}(1-x)^{s}} d x
\end{aligned}
\end{equation}

Performing this procedure for $a_{p} \ldots a_{0}$, summing the $\mathrm{p}+1$ terms, and factoring out $\frac{(a+1)_{k n}(b+1)_{s n}}{(a+b+2)_{(k+s) n}(z)^{n}}$ gives equation (1.4), which completes the proof.

One final pochhammer identity that will be encountered frequently in this paper is:

$$
(a)_{n k}=\left(\prod_{y=0}^{k-1}\left(\frac{a+y}{k}\right)_{n}\right) k^{k n}
$$

Another formula that will be used extensively later is the Euler integral:

\begin{equation}
2 F 1\left\{x_{1}, x_{2} ; y_{1} ; z\right\}=\frac{\Gamma\left(y_{1}\right)}{\Gamma\left(x_{2}\right) \Gamma\left(y_{1}-x_{2}\right)} \int_{0}^{1} \frac{x^{x_{2}-1}(1-x)^{y_{1}-x_{2}-1}}{(1-z x)^{x_{1}}} d x
\end{equation}

\subsection{Deriving Formulas for Pi}
We shall use (1.4) to prove (1.1). Consider a hypergeometric function for arcsine converted into integral form via (2.6):

\begin{equation}
2 F 1\left\{1, \frac{1}{2} ; \frac{3}{2} ; \frac{w^{2}}{w^{2}-1}\right\}=\left(\frac{\sin ^{-1}(w)}{w}\right) \sqrt{\left(1-w^{2}\right)}=\frac{1}{2} \int_{0}^{1} \frac{x^{-1 / 2}}{1-\left(\frac{x \cdot w^{2}}{w^{2}-1}\right)} d x
\end{equation}

Then let $\mathrm{w}=1 / 2$ so that $\frac{\pi \sqrt{3}}{3}=\int_{0}^{1} \frac{x^{-1 / 2}}{1+\frac{x}{3}} d x$

Using (1.4) let $\mathrm{k}=1$ and $\mathrm{s}=2$. Then let $\mathrm{a}=-1 / 2 \mathrm{~b}=0$ and $P(x)=1+\frac{x}{3}$ and find $\mathrm{Q}(\mathrm{x})$ such that:

$$
\int_{0}^{1} \frac{x^{-1 / 2} Q(x)}{z-x(1-x)^{2}} d x=\int_{0}^{1} \frac{x^{-1 / 2}}{1+\frac{x}{3}} d x
$$

Expand $z-x(1-x)^{2}$ and perform polynomial division $\frac { x } { 3 } + 1 \longdiv { - x ^ { 3 } + 2 x ^ { 2 } - x + z }$ A value of $\mathrm{z}=-48$ and $Q(x)=-3 x^{2}+15 x-48$ is obtained. Thus, from which (1.1) follows, completing the proof:

\begin{equation}
\int_{0}^{1} \frac{x^{-1 / 2}\left(x^{2}-5 x+16\right)}{-48-x(1-x)^{2}} d x=\frac{-1}{3} \int_{0}^{1} \frac{x^{-1 / 2}}{1+\frac{x}{3}} d x=\frac{-\pi \sqrt{3}}{9}
\end{equation}

(2.7) can be generalized, in which  (1.1) immediately follows for $w=1/2$ :

\begin{equation}
\begin{aligned}
&  \sum_{n=0}^{\infty}\left(\frac{4w^6}{27(w^2-1)}\right)^n\frac{(1)_{n}(1/2)_{n}\left(n(4w^4+6w^2-18)+2w^4+5w^2-15\right)}{(11/6)_{n}(7/6)_{n}} = \\ 
& \left(-\frac{15\sin ^{-1}(w)}{w}\right) \sqrt{\left(1-w^{2}\right)}
\end{aligned}
\end{equation}

This follows from:

\begin{equation}
 \frac{w^2}{1-w^2x} = \frac{x^2+x(1/w^2-1)+(1-w^2)/w^4}{(w^2-1)/w^6-x^2(1-x)} 
\end{equation}

$$$$

Equation (1.1) converges at a rate of $\log (324)$ or 2.5 digits per term using an arcsine formula for $1 / 2$, which to the best  of my knowledge is faster than any previously known formula of this type bi-3F2 series. The result is a rather compact formula that converges much faster than similar formula in the literature.  Part of the reason it converges so fast is because two of the pochammer symbols  cancel out for (2.5) and fractions cancel out for (2.4) owing to the  factorization of 6 into two unique primes (e.g. $6n+3=2n+1$), so instead of it being a 4F3 series, it's a 3F2 series yet retains the accelerated convergence rate. 

Like similar binomial  $\pi$ formulas, aesthetically it bears a resemblance to the Ramanujan $\pi$ formulas, but is derived through elementary means instead of elliptic integrals, and not the reciprocal of $\pi$. Like the Ramanujan formulas, (2.9) can be split into the sum of two 3F2 hypergeometric series. 

I conjecture 1.1 is the fastest $\pi$ formula that can be expressed  in the form $c \pi= \sum_{n=0}^{\infty} Q(n) (a+bn)z^{n} $ in which $c$ is algebraic and $a,b,z$ are rational, and Q(n) is some array of pochammer symbols for a 3F2 hypergeoemtric series.  

A faster formula can be obtained that adds five digits per term. Let $s=4, k=2, a=-1 / 2$, and $\mathrm{b}=0 . P(x)=\left(1+\frac{x}{3}\right)$ Then expand $z-x^{2}(1-x)^{4}$ and performing division for $\mathrm{Q}(\mathrm{x})$ the following identity is obtained:

\begin{equation}
\int_{0}^{1} \frac{x^{-1 / 2}\left(-x^{5}+7 x^{4}-27 x^{3}+85 x^{2}-256 x+768\right)}{2304-x^{2}(1-x)^{4}} d x=\frac{\pi \sqrt{3}}{9}
\end{equation}

After some labor, we have:

\begin{equation}
\pi=\frac{\sqrt{3}}{6^{5}} \sum_{n=0}^{\infty}\left(\frac{[(4 n) !]^{2}(6 n) !}{(2 n) !(12 n) !(9)^{n+1}}\right)\left(\frac{127169}{12 n+1}-\frac{1070}{12 n+5}-\frac{131}{12 n+7}+\frac{2}{12 n+11}\right)
\end{equation}
 
To derive another pi formula let $\mathrm{s}=2, \mathrm{k}=2, \mathrm{a}=-1 / 2, \mathrm{~b}=0$ and $P(x)=1+\frac{x}{3}$ to obtain:

\begin{equation}
\pi 2^{10} \sqrt{3}=\sum_{k=0}^{\infty} \frac{1}{9^k\left(\begin{array}{c}8 k \\ 4 k\end{array}\right)}\left(\frac{5717}{8 k+1}-\frac{413}{8 k+3}-\frac{45}{8 k+5}+\frac{5}{8 k+7}\right)
\end{equation}

This is similar to [4]

\section{Multivariable Values of P(x) and Q(x) for Infinite Series}
\numberwithin{equation}{section}

Binomial identities can be found using equations (1.3) and (1.4), but extending it to multiple variables. In this section $\mathrm{Q}(\mathrm{x})$ becomes $\mathrm{Q}(\mathrm{x}, \mathrm{w}), \mathrm{P}(\mathrm{x})$ becomes $\mathrm{P}(\mathrm{x}, \mathrm{w})$, and $\mathrm{z}$ becomes a function of $\mathrm{w}$.

To begin, using equation (1.3) let $\mathrm{a}=0, \mathrm{~b}=0, \mathrm{k}=1, \mathrm{~s}=1, \mathrm{z}=\mathrm{w}$ to obtain the identity:

\begin{equation}
\int_{0}^{1} \frac{1}{w-x(1-x)} d x=\sum_{n=0}^{\infty} \frac{1}{w^{n+1}\left(\begin{array}{l}
2 n \\
n
\end{array}\right)(2 n+1)}
\end{equation}

By evaluating the left hand integral we see that (3.1) is equal to [3]:

\begin{equation}
\frac{4 \tan ^{-1}(1 / \sqrt{4 w-1})}{\sqrt{4 w-1}}
\end{equation}

Do derive a faster series for (3.1) let $\mathrm{k}=3, \mathrm{~s}=3, \mathrm{a}=0, \mathrm{~b}=0, \mathrm{z}=\mathrm{z}$, and expand $\mathrm{z}-x^{3}(1-x)^{3}$

And let $\mathrm{P}(\mathrm{x}, \mathrm{w})=x^{2}-x+w$. Performing polynomial division, we need to find $\mathrm{Q}(\mathrm{x}, \mathrm{w})$ and a new value of $\mathrm{z}$ in terms of $\mathrm { w } \cdot ( x ^ { 2 } - x + w ) \longdiv { ( x ^ { 6 } - 3 x ^ { 5 } + 3 x ^ { 4 } - x ^ { 3 } + z ) }$

We find: $Q(x, w)=\left(x^{4}-2 x^{3}+(1-w) x^{2}+x w+w^{2}\right), z=w^{3}$ Therefore:

\begin{equation}
\int_{0}^{1} \frac{\left(x^{4}-2 x^{3}+(1-w) x^{2}+x w+w^{2}\right)}{w^{3}-x^{3}(1-x)^{3}} d x=\int_{0}^{1} \frac{1}{w-x(1-x)} d x
\end{equation}

With regard to equation (3.1) we can let $w=w^{-2}$,divide by $\mathrm{w}$, and take the derivative in terms of w such that:

\begin{equation}
\sum_{n=0}^{\infty} \frac{w^{2 n+1}}{\left(\begin{array}{l}
2 n \\
n
\end{array}\right)(2 n+1)} \frac{d}{d w}=\sum_{n=0}^{\infty} \frac{w^{2 n}}{\left(\begin{array}{l}
2 n \\
n
\end{array}\right)}=\frac{w^{2}}{4-w^{2}}+\frac{4 w \tan ^{-1}\left(w / \sqrt{4-w^{2}}\right)}{\left(4-w^{2}\right)^{3 / 2}}
\end{equation}

An arbitrary number of integrals and derivatives of (3.4) can be taken in terms of $\mathrm{w}$ to obtain various central binomial identities.

Applying formula (1.4) on (3.3), and after some labor obtain:

(The variable $r$ denotes the powers of $n$ that results from taking multiple derivatives)

\begin{equation}
\begin{aligned}
& 4 \sum_{n=1}^{\infty} \frac{w^{n}}{\left(\begin{array}{l}
2 n \\
n
\end{array}\right) n^{r}}= \\
& 2 w^{2}+\frac{2 w^{2}}{3 \cdot 2^{r}}+\sum_{n=1}^{\infty} \frac{w^{3 n}}{\left(\begin{array}{l}
6 n \\
3 n
\end{array}\right)}\left(\frac{w^{2}\left(9 n^{2}+9 n+2\right)}{(3 n+2)^{r}(6 n+3)(6 n+1)}+\frac{2 w(3 n+1)}{(6 n+1)(3 n+1)^{r}}+\frac{4}{(3 n)^{r}}\right)
\end{aligned}
\end{equation}

Let $\mathrm{r}=0$ and $w=w^{2}$ in (3.5). Then divide (3.5) by $\mathrm{w}$ and take multiple derivatives such that the denominators vanish, ensuring you divide (3.5) by $\mathrm{w}$ each time before the derivative is taken. Then finally let $\mathrm{w}=1$ to obtain the identity:

\begin{equation}
3 \sum_{n=1}^{\infty} \frac{63 n^{2}-27 n+4}{\left(\begin{array}{l}
6 n \\
3 n
\end{array}\right)}=16 \sum_{n=1}^{\infty} \frac{n^{2}}{\left(\begin{array}{l}
2 n \\
n
\end{array}\right)}-32 \sum_{n=1}^{\infty} \frac{n}{\left(\begin{array}{l}
2 n \\
n
\end{array}\right)}+12 \sum_{n=1}^{\infty} \frac{1}{\left(\begin{array}{l}
2 n \\
n
\end{array}\right)}
\end{equation}

Using (3.4) and taking multiple derivatives we arrive at:

\begin{equation}
3 \sum_{n=1}^{\infty} \frac{63 n^{2}-27 n+4}{\left(\begin{array}{l}
6 n \\
3 n
\end{array}\right)}=\frac{40 \pi \sqrt{3}}{81}+4
\end{equation}

A faster converting series can be obtained by letting $\mathrm{k}=5, \mathrm{~s}=5$ such that

\begin{equation}
\begin{aligned}
& \int_{0}^{1} \frac{Q(x, w)}{w^{5}-x^{5}(1-x)^{5}} d x=\int_{0}^{1} \frac{1}{P(x, w)} d x \\
& P(x, w)=w-x(1-x) \\
& Q(x, w)=\left(\begin{array}{l}
x^{8}-4 x^{7}+(6-w) x^{6}+(3 w-4) x^{5}+\left(w^{2}-3 w+1\right) x^{4} \\
+\left(w-2 w^{2}\right) x^{3}+\left(w^{2}-w^{3}\right) x^{2}+w^{3} x+w^{4}
\end{array}\right)
\end{aligned}
\end{equation}

Using formula (1.4) and after some labor the identity is obtained:

\begin{equation}
\begin{aligned}
& \sum_{n=0}^{\infty} \frac{213125 n^{4}-278000 n^{3}+139975 n^{2}-26800 n+1596}{\left(\begin{array}{l}
10 n \\
5 n
\end{array}\right)}= \\
& 16 \sum_{n=0}^{\infty} \frac{16 n^{4}-128 n^{3}+344 n^{2}-352 n+105}{\left(\begin{array}{l}
2 n \\
n
\end{array}\right)}=\frac{1120 \pi \sqrt{3}}{81}+1728
\end{aligned}
\end{equation}

\section{Hypergeometric Function Transformations for the Gamma Function}
\numberwithin{equation}{section}

Infinite series for the gamma function such as (1.2) can be derived using identities (1.3) and (1.4). In this section we'll also prove a
quadratic hypergeometric transformation that is a byproduct of this method, as well as other identities.

Begin with Kummer’s formula [5] and express it in integral form via (2.6):

\begin{equation}
2 F 1\{1, h ; 2-h ;-1\}=\frac{\Gamma(2-h) \Gamma(3 / 2)}{\Gamma(3 / 2-h) \Gamma(2)}=\frac{\Gamma(2-h)}{\Gamma(2-2 h) \Gamma(h)} \int_{0}^{1} \frac{x^{h-1}(1-x)^{1-2 h}}{1+x} d x
\end{equation}

Since $-2-x(1-x)=(x-2)(x+1)$ we can let $\mathrm{k}=1, \mathrm{~s}=1, \mathrm{z}=-2$ and $Q(x)=x-2$ And after isolating the integral in (4.1) we have the identity:

\begin{equation}
\int_{0}^{1} \frac{x^{h-1}(1-x)^{1-2 h}(x-2)}{-2-x(1-x)} d x=\int_{0}^{1} \frac{x^{h-1}(1-x)^{1-2 h}}{x+1} d x=\frac{\Gamma(2-2 h) \Gamma(h) \sqrt{\pi}}{2 \Gamma(3 / 2-h)}
\end{equation}

The leftmost integral in (4.3) can be evaluated using (1.4). setting $a_{0}=-2 a_{1}=1$ $\mathrm{a}=\mathrm{h}-1 \mathrm{~b}=1-2 \mathrm{~h}$ :

$$
\frac{\Gamma(h) \Gamma(2-2 h)}{-2 \Gamma(2-h)} \sum_{n=0}^{\infty} \frac{(h)_{n}(2-2 h)_{n}}{(2-h)_{2 n}(-2)^{n}}\left(\frac{h+n}{2-h+2 n}-2\right)=\frac{\sqrt{\pi} \Gamma(h) \Gamma(2-2 h)}{2 \Gamma(3 / 2-h)}
$$

Using (2.5), we obtain the series computing values of the gamma function:

\begin{equation}
\frac{h(h+1) 2^{2 h-1}[\Gamma(h)]^{2}}{\Gamma(2 h)}=\sum_{n=0}^{\infty} \frac{(1-h)_{n}(2 h)_{n}(3 n+3 h+1)}{(-8)^{n}(h / 2+1)_{n}(h / 2+3 / 2)_{n}}
\end{equation}

To derive (1.2) let $\mathrm{h}=1 / 3$ in equation (4.3) and simplify:

$$
\frac{3^{1 / 2} 2^{2 / 3}[\Gamma(1 / 3)]^{3}}{9 \pi}=\sum_{n=0}^{\infty} \frac{(2 / 3)_{n}(2 / 3)_{n}(3 n+2)}{(-8)^{n}(7 / 6)_{n}(5 / 3)_{n}}
$$

But $(2 / 3)_{n}(3 n+2)=2(5 / 3)_{n}$ so we can write:

\begin{equation}
\frac{3^{1 / 2} 2^{2 / 3}[\Gamma(1 / 3)]^{3}}{18 \pi}=\sum_{n=0}^{\infty} \frac{(2 / 3)_{n}}{(-8)^{n}(7 / 6)_{n}}=2 F 1\left\{1, \frac{2}{3} ; \frac{7}{6} ; \frac{-1}{8}\right\}
\end{equation}

(4.4)Is a quadratically transformed hypergeometric function evaluated at a specific point. This to the best of my knowledge is a new derivation of this $\Gamma(1 / 3)$ series, which does not involve either manipulation of the hypergeometric function, transformation of the Euler integral, or evaluating the hypergeoemtric differential equation itself, but rather by applying the Beta method to Kummer’s formula . 

Via (2.6):

$$
\frac{3^{1 / 2} 2^{2 / 3}[\Gamma(1 / 3)]^{3}}{18 \pi}=\frac{\Gamma(7 / 6)}{\Gamma(1 / 2) \cdot \Gamma(2 / 3)} \int_{0}^{1} \frac{x^{-1 / 3}(1-x)^{-1 / 2}}{1+x / 8} d x
$$

We obtain the simplification:

\begin{equation}
\frac{4 \cdot 3^{1 / 2} \cdot \pi}{9}=\int_{0}^{1} \frac{x^{-1 / 3}(1-x)^{-1 / 2}}{1+x / 8} d x
\end{equation}

Letting

$$
\begin{aligned}
& k=1, s=1, a=-1 / 3, b=-1 / 2, P(x)=(1+x / 8) \\
& Q(x)=(8 x-72), z=-72, a_{1}=8, a_{0}=-72
\end{aligned}
$$

We arrive at the identity:

$$
\int_{0}^{1} \frac{x^{-1 / 3}(1-x)^{-1 / 2}}{1+x / 8} d x=\int_{0}^{1} \frac{x^{-1 / 3}(1-x)^{-1 / 2}(8 x-72)}{-72-x(1-x)} d x
$$

And  obtain   (1.2).
 
\section{Higher-Order Constants}
\numberwithin{equation}{section}

There series we've been deriving so far involve constants that can be expressed though the hypergeometric function of the form $2 F 1\{1, x ; y ; z\}$. Such constants involve logarithms, root extractions, the gamma function, and inverse trigonometric functions. Higher order constants such as Catalans's Constant, Zeta 3, etc can't be expressed through $2 F 1\{1, x ; y ; z\}$, but require a more generalized hypergeometric function of the form: ${ }_{q+1} F_{q}\left\{1, x_{1}, \ldots, x_{q} ; y_{1} ; \ldots ; y_{q} ; z\right\}$.

Unfortunately, the algorithm (1.4) isn't as efficient at generating rapidly converging series for constants that are only expressible though a hypergeometric function where $q>1$. To the best of my knowledge for (1.4) to work on higher order constants such as Zeta 3, the value of ' $s$ ' must be set to zero. Accelerating the convergence of a q>1 hypergeometric function via (1.4) by setting $\mathrm{s}=0$ equates to taking multiple terms of

${ }_{q+1} F_{q}\left\{1, x_{1}, \ldots, x_{q} ; y_{1} ; \ldots ; y_{q} ; z\right\}$ ' $\mathrm{k}$ ' terms at a time. So if $\mathrm{k}=2$, the series generated by (1.4) would converge at a rate of $z^{2}$.

In this next example we will show how to accelerate the convergence of a hypergeometric function of the form $3 F 2\left\{1, x_{1}, x_{2} ; y_{1} ; y_{2} ; z\right\}$ using (1.4), and then apply it to derive faster converging infinite series for Catalan's Constant and $\Gamma(1 / 4)$.

Begin with the integrals:
\begin{equation}
\begin{aligned}
 3 F 2\left\{1, x_{1}, x_{2} ; y_{1} ; y_{2} ; z\right\}=\\ 
 \frac{\Gamma\left(y_{2}\right)}{\Gamma\left(x_{2}\right) \Gamma\left(y_{2}-x_{2}\right)} \int_{0}^{1} x^{x_{2}-1}(1-x)^{y_{2}-x_{2}-1} 2 F 1\left\{1, x_{1} ; y_{1} ; x z\right\} d x 
\end{aligned}
\end{equation}

\begin{equation}
2 F 1\left\{1, x_{1} ; y_{1} ; x z\right\}=\frac{\Gamma\left(y_{1}\right)}{\Gamma\left(x_{1}\right) \Gamma\left(y_{1}-x_{1}\right)} \int_{0}^{1} \frac{y^{x_{1}-1}(1-y)^{y_{1}-x_{1}-1}}{1-y x z} d y
\end{equation}

Using the integral (5.2), choosing the values below for $\mathrm{k}$ and $\mathrm{s}$, and performing polynomial long division we plug the following into (1.4):

 $$
\begin{aligned}
& s=0, k=2, b=\left(y_{1}-x_{1}-1\right), a=\left(x_{1}-1\right), z=\frac{1}{x^{2} z^{2}} \\
& P(x, y, z)=(1-x y z), Q(x, y, z)=\left(\frac{y}{x z}+\frac{1}{x^{2} z^{2}}\right)
\end{aligned}
 $$

And obtain:

\begin{equation}
\begin{aligned}
& \int_{0}^{1} \frac{y^{x_{1}-1}(1-y)^{y_{1}-x_{1}-1}}{P(x, y, z)} d y=\int_{0}^{1} \frac{y^{x_{1}-1}(1-y)^{y_{1}-x_{1}-1} Q(x, y, z)}{\frac{1}{x^{2} z^{2}}-y^{2}} d y= \\
& \frac{x^{2} z^{2} \cdot \Gamma\left(x_{1}\right) \Gamma\left(y_{1}-x_{1}\right)}{\Gamma\left(y_{1}\right)} \sum_{n=0}^{\infty} \frac{\left(x_{1}\right)_{2 n}(x z)^{2 n}}{\left(y_{1}\right)_{2 n}}\left(\frac{x_{1}+2 n}{x z\left(y_{1}+2 n\right)}+\frac{1}{x^{2} z^{2}}\right)
\end{aligned}
\end{equation}

Plugging (5.3) into (5.2), and then plugging that result into (5.1) and simplifying we get:

\begin{equation}
\begin{aligned}
& \frac{\Gamma\left(y_{2}\right)}{\Gamma\left(x_{2}\right) \Gamma\left(y_{2}-x_{2}\right)} \int_{0}^{1} x^{x_{2}-1}(1-x)^{y_{2}-x_{2}-1} \sum_{n=0}^{\infty} \frac{\left(x_{1}\right)_{2 n}(x z)^{2 n}}{\left(y_{1}\right)_{2 n}}\left(\frac{x z\left(x_{1}+2 n\right)}{y_{1}+2 n}+1\right) d x \\
& =3 F 2\left\{1, x_{1}, x_{2} ; y_{1} ; y_{2} ; z\right\}
\end{aligned}
\end{equation}

Then we extract and convert the following integrals from (5.4) into pochhammer \& gamma notation via the beta function:

\begin{equation}
\begin{aligned}
& \int_{0}^{1} x^{2 n+x_{2}-1}(1-x)^{y_{2}-x_{2}-1} d x=\frac{\Gamma\left(2 n+x_{2}\right) \Gamma\left(y_{2}-x_{2}\right)}{\Gamma\left(2 n+y_{2}\right)}=\frac{\left(x_{2}\right)_{2 n} \Gamma\left(x_{2}\right) \Gamma\left(y_{2}-x_{2}\right)}{\left(y_{2}\right)_{2 n} \Gamma\left(y_{2}\right)} 
\end{aligned}
\end{equation}

Then plugging (5.5) back into (5.4) and simplifying gives the final result:

\begin{equation}
3 F 2\left\{1, x_{1}, x_{2} ; y_{1} ; y_{2} ; z\right\}=\sum_{n=0}^{\infty} \frac{\left(x_{1}\right)_{2 n}\left(x_{2}\right)_{2 n}(z)^{2 n}}{\left(y_{1}\right)_{2 n}\left(y_{2}\right)_{2 n}}\left(\frac{z\left(x_{1}+2 n\right)\left(x_{2}+2 n\right)}{\left(y_{1}+2 n\right)\left(y_{2}+2 n\right)}+1\right)
\end{equation}

The procedure used to derive (5.6) can be generalized:

\begin{equation}
{ }_{q+1} F_{q}\left\{1, x_{1}, \ldots, x_{q} ; y_{1} ; \ldots ; y_{q} ; z\right\}=\sum_{n=1}^{\infty}\left(\left(\left(z \prod_{g=1}^{q} \frac{x_{g}+2 n}{y_{g}+2 n}\right)+1\right) \prod_{g=1}^{q} \frac{\left(x_{g}\right)_{2 n}}{\left(y_{g}\right)_{2 n}}\right)
\end{equation}

With regards to Catalan's Constant it is well-known that

\begin{equation}
3 F 2\left\{1,1, \frac{1}{2} ; \frac{3}{2} ; \frac{3}{2} ; \frac{1}{4}\right\}=\sum_{n=0}^{\infty} \frac{1}{(2 n+1)^{2}\left(\begin{array}{l}
2 n \\
n
\end{array}\right)}=\frac{\pi}{3} \ln (2-\sqrt{3})+\frac{8}{3} \sum_{n=0}^{\infty} \frac{(-1)^{n}}{(2 n+1)^{2}}
\end{equation}

Plugging the hyqpergeometric terms of (5.8) into (5.6) and simplifying gives an infinite series for (5.8) that converges twice as fast:

\begin{equation}
\sum_{n=0}^{\infty} \frac{1}{(2 n+1)^{2}\left(\begin{array}{l}
2 n \\
n
\end{array}\right)}=\sum_{n=0}^{\infty} \frac{40 n^{2}+54 n+19}{((4 n+1)(4 n+3))^{2}\left(\begin{array}{l}
4 n \\
2 n
\end{array}\right)}
\end{equation}

Using a transformation on (5.2) letting $\mathrm{k}=3$ gives:

\begin{equation}
4 \sum_{n=0}^{\infty} \frac{1}{(2 n+1)^{2}\left(\begin{array}{l}
2 n \\
n
\end{array}\right)}=\sum_{n=0}^{\infty} \frac{6804 n^{4}+17172 n^{3}+15903 n^{2}+6405 n+956}{((6 n+1)(6 n+3)(6 n+5))^{2}\left(\begin{array}{l}
6 n \\
3 n
\end{array}\right)}
\end{equation}

In addition, formulas for $\Gamma(1 / 4)$ can be derived using (5.6) based on accelerating a quadratic transformation. Begin with the following quadratically transformed hypergeometric functions:

\begin{equation}
2 F 1\left\{\frac{1}{4}, \frac{1}{4} ; 1 ; \frac{-1}{8}\right\}=\frac{\sqrt{\pi}}{2^{1 / 4}[\Gamma(3 / 4)]^{2}} 
\end{equation}
\begin{equation}
2 F 1\left\{\frac{1}{4}, \frac{3}{4} ; 1 ; \frac{1}{9}\right\}=\frac{\sqrt{3 \pi}[\Gamma(1 / 4)]^{2}}{4 \pi^{2}}
\end{equation}

(Note: Unlike equation (4.4), (5.11) and (5.12) can’t be derived with (1.4). Proofs regarding quadratic hypergeometric transformations is given in [7].)

However, (5.11) is equal to $3 F 2\left\{1, \frac{1}{4}, \frac{1}{4} ; 1 ; 1 ; \frac{-1}{8}\right\}$ and (5.12) can be written as $3 F 2\left\{1, \frac{1}{4}, \frac{3}{4} ; 1 ; 1 ; \frac{1}{9}\right\}$. Plugging these values into (5.6) and simplifying via (2.5) gives

\begin{equation}
\frac{36 \sqrt{3}[\Gamma(1 / 4)]^{2}}{\pi^{3 / 2}}=\sum_{n=0}^{\infty} \frac{\left(640 n^{2}+608 n+147\right)(8 n) !}{[(2 n+1) !]^{2}(4 n) ! 24^{4 n}}
\end{equation}

\begin{equation}
\frac{128 \pi^{1 / 2}}{2^{1 / 4}[\Gamma(3 / 4)]^{2}}=\sum_{n=0}^{\infty}\left(448 n^{2}+448 n+127\right)\left(\frac{(1 / 8)_{n}(5 / 8)_{n}}{2^{n}(2 n+1) !}\right)^{2}
\end{equation}

\section*{References:}
[1] Boris Gour'evitch, Jes'us Guillera Goyanes, Construction Of Binomial Sums For pi And Polylogarithmic Constants Inspired By BBP Formulas (2006).

[2] D. Bailey, P. Borwein and S. Plouffe, On The Rapid Computation of Various Polylogarithmic Constants, Math. Comp., 66(1997), 903-913.

[3] Travis Sherman, Summation of Glaisher- and Apéry-like Series, University of Arizona (2000),\href{http://math.arizona.edu/}{http://math.arizona.edu/} ura/001/sherman.travis/series.pdf

[4] G. Almkvist, C. Krattenthaler and J. Petersson, Some new formulas for \_, preprint, Matematiska Institutionen, Lunds Universitet, Sweden. and Institut fur Mathematik der Universitat Wien, Austria (2001), \newline \href{http://www.arxiv.org/pdf/math.NT/0110238}{http://www.arxiv.org/pdf/math.NT/0110238}

[5] "List of hypergeometric identities" Wikipedia, The Free Encyclopedia. Wikimedia Foundation, Inc. 10 Aug. 2004.

\href{http://en.wikipedia.org/wiki/List_of_hypergeometric_identities}{http://en.wikipedia.org/wiki/List\_of\_hypergeometric\_identities}

[6] Fabrice Bellard, \_ page, \href{http://fabrice.bellard.free.fr/pi/}{http://fabrice.bellard.free.fr/pi/}.

[7] Earl D. Rainville, Ph.D. "Special Functions", page 65, The MacMillan Company, NY., 1960.

[8] Hakimoglu, Cetin, "An Algorithm for the Derivation of Rapidly Converging Infinite Series for Universal Mathematical Constants" (2009,2021). Available at SSRN: \href{https://ssrn.com/abstract=3919892}{https://ssrn.com/abstract=3919892}  

[9] John M. Campbell et al., "On Some Series Involving the Binomial
Coefficients", 2019, https://arxiv.org/pdf/2306.16889.pdf

[10]  Wenchang Chu, "DOUGALL’S BILATERAL 2H2-SERIES
AND RAMANUJAN-LIKE FORMULAE", 2011, \newline https://www.ams.org/journals/mcom/2011-80-276/S0025-5718-2011-02474-9/S0025-5718-2011-02474-9.pdf

[11] Wenchang Chu, "Ramanujan-like formulae for and via Gould–Hsu inverse series relations", 2011, https://link.springer.com/article/10.1007/s11139-020-00337-z 

[12] Zhi-Wei Sun, "New series involving binomial coefficients", \newline 2023, \href{https://arxiv.org/abs/2307.03086}{https://arxiv.org/abs/2307.03086}

[13] K. A. DRIVER, "AN INTEGRAL REPRESENTATION OF SOME HYPERGEOMETRIC
FUNCTIONS", 2006, \newline https://www.emis.de/journals/ETNA/vol.25.2006/pp115-120.dir/pp115-120.pdf 

[14] Necdet Batir, "On the series", 2005, https://arxiv.org/pdf/math/0512310.pdf

\end{document}